\documentclass[10pt,a4paper]{article}
\usepackage[utf8]{inputenc}
\usepackage[english]{babel}
\usepackage{amsmath}
\usepackage{amsfonts}
\usepackage{amssymb}
\usepackage{amsthm}
\usepackage{lipsum}
\usepackage[left=2cm,right=2cm,top=2cm,bottom=2cm]{geometry}

\newcommand*\lap{\Delta}
\newcommand{\Sp}{\mathbb{S}}

\newcommand{\R}{\mathbb{R}}
\newcommand{\eqdef}{\overset{def}{=}}	
\newcommand{\vol}{\mathrm{vol}}

\newtheorem{rem}{Remark}

\newtheorem{thm}{Theorem}

\title{Integrals of spherical harmonics with Fourier exponents in multidimensions}
\author{F. O. Goncharov\thanks{CMAP, Ecole Polytechnique, CNRS, Universit\'{e} Paris-Saclay, 91128, Palaiseau, France; email: fedor.goncharov.ol@gmail.com}}
\begin{document}
\maketitle
\abstract{We consider integrals of spherical harmonics with Fourier exponents on the sphere $\Sp^{n}, \, n\geq 1$. Such transforms arise in the framework of the theory of weighted Radon transforms and vector diffraction in electromagnetic fields theory. We give analytic formulas for these integrals, which are exact up to multiplicative constants. These constants depend on choice of basis on the sphere. In addition, we find these constants explicitly for the class of harmonics arising in the framework of the theory of weighted Radon transforms. We also suggest formulas for finding these constants for the general case.}

\section{Introduction}
We consider the integrals
\begin{equation}\label{integral_def}
	I_{k}^{m}(p,\rho) \eqdef \int\limits_{\Sp^{n}} Y_{k}^{m}(\theta)e^{i\rho(p,\theta)}\, d\theta,
	\, p\in \Sp^{n}, \, \rho \geq 0, \, n\geq 1,
\end{equation}
where $\{Y_{k}^{m} \, | \,  k\in \mathbb{N}\cup \{0\}, \, m = \overline{1, a_{k,n+1}}\}$ is orthonormal basis of spherical harmonics  on $\Sp^{n}\subset \mathbb{R}^{n+1}$ (see e.g. \cite{morimoto1998analytic}, \cite{stein2016introduction}), $a_{k,n}$ is defined  as follows:
\begin{equation}\label{eigensp_dim}
	a_{k,n+1} = {{n+k}\choose{k}} + {{n+k-2}\choose{k-2}}, \, a_{0,n} = 1,\, a_{1,n} = n,
\end{equation}
where 
\begin{equation}
	{{n}\choose{k}} = \dfrac{n!}{k!(n-k!)}, \, n,k\in\mathbb{N}\cup \{0\}.
\end{equation}
\par We recall that spherical harmonics $\{Y^m_k\}$ are eigenfunctions of the spherical Laplacian $\lap_{\Sp^n}$ and 
the following identity holds (see e.g. \cite{morimoto1998analytic}, \cite{stein2016introduction}):
\begin{equation}\label{eigen_sp_lap}
	\lap_{\Sp^n}Y_{k}^{m} = -k(n+k-1)Y_k^{m}, \, m = \overline{1,a_{k,n+1}},
\end{equation}
where $a_{k,n}$ is defined in \eqref{eigensp_dim}.
\par Integrals $I_k^m$ arise, in particular, in connection with iterative inversions of the weighted Radon transforms in dimension $d=n+1=2$; see \cite{kunyansky1992generalized}. In addition, 
an exact (up to multiplicative coefficient depending on $k$) analytic formula for \eqref{integral_def} was given in \cite{kunyansky1992generalized} for $d=n+1=2, \, k=2j, \, j\in \mathbb{N}\cup \{0\}$. 
\par We recall that weighted Radon transform operator $R_W$ is defined (in dimension $d=n+1$) as follows (see e.g.
 \cite{goncharov2016analog},  \cite{kunyansky1992generalized}):
\begin{equation}\label{radon_def}
	R_{W}f(s,\theta) \eqdef \int\limits_{x\theta =s }W(x,\theta) f(x)\, dx, \, 
	s\in \R, \, \theta\in \Sp^n,
\end{equation}
where $W$ is the weight function on $\R^{n+1}\times \Sp^{n}$, $f$ is a test-function.
\par The present work is strongly motivated by the fact that integrals $I_k^m$ also arise in the theory of weighted Radon transforms defined by \eqref{radon_def} for higher dimensions $d=n+1\geq 3$. This issue will be presented in detail in the subsequent work \cite{gon2016nov}. 
\par On the other hand, in \cite{neves2006analytical} integrals $I_k^m$ were considered for the case of $n=2$  in connection with vector diffraction in electromagnetic theory and exact analytic formulas were given for this case. 
\par In addition, for the case of dimension $n=2$ more general forms of integrals $I_k^m$  were considered in the recent work \cite{adkins2013three}. In particular, the results of \cite{adkins2013three} coincide with the results of the present work for the case of dimension $n=2$. 

\par In the present work we prove that 
\begin{equation}\label{first_result}
I_k^m(p,\rho)= c(m,k,n) Y_k^m(p) \rho^{(1-n)/2}J_{k + \frac{n-1}{2}}(\rho),
\end{equation}
where $J_r(\cdot)$ is the $r$-th Bessel function of the first kind, $c(m,k,n)$  is a constant which depends on indexes $m,k$ of spherical harmonic $Y_k^m$ and on dimension $n$; see Theorem 1 in Section~\ref{sect_main_res}.
\par This result is new for the case of $d=n+1=2$ for odd $k$ and for $d=n+1>3$ in general.
\par In the framework of applications to the theory of weighted Radon transforms, integrals $I_k^m$ arise for the case of even $k=2j, \, j\in \mathbb{N}\cup \{0\}$; see formula (7) of  \cite{kunyansky1992generalized} for  $n=1$ and subsequent work \cite{gon2016nov} for $n\geq 2$. 
For $k=2j$ we find explicitly the constants $c(m,2j,n)$ arising in \eqref{first_result}; see Theorem~\ref{thm_target} in Section~\ref{sect_main_res}. 
\par It is interesting to note that the constants $c(m,2j,n)$ are expressed via the eigenvalues of the Minkowski-Funk transform $\mathcal{M}$ on $\Sp^n$, where operator $\mathcal{M}$ is defined as follows (see e.g. \cite{funk1913}, \cite{gel_fand2003selected}):  
\begin{equation}\label{minkowski_funk}
	\mathcal{M}[f](p) = \int\limits_{\Sp^n, \, (\theta p) = 0}f(\theta)\, d\theta, \, 
	p\in \Sp^n,
\end{equation}
where $f$ is an even test-function on $\Sp^n$; see Section~\ref{sect_main_res} for details.
\par In Section~\ref{sect_proofs} we give proofs of  Theorem~\ref{thm_target} and Remark~\ref{lem_general_sum}.

\section{Main results}\label{sect_main_res}

\begin{thm}\label{thm_target}  Let $I_{k}^{i}(p,\rho)$ be defined by \eqref{integral_def}. Then:
\par (i) The following formulas hold:
\begin{align}\label{integral_claims}
	&I_k^{m}(p,\rho) = c(m,k,n) Y_{k}^{m}(p) \rho^{(1-n)/2}J_{k+ \frac{n-1}{2}}(\rho),\\
	\label{integral_sign_claims}
	&I_k^m(p, -\rho) = (-1)^k I_{k}^m(p,\rho),\\ \nonumber
	&p\in \Sp^n, \, \rho \in \mathbb{R}_+ = [0, +\infty),
\end{align}
where $J_r(\rho)$ is the standard $r$-th Bessel function of the first kind, $c(m, k, n)$ depends only on integers $m,\,k, \,n$ for fixed orthonormal basis $\{Y^m_k\}$.
\par (ii) In addition, for $k = 2j, \, j \in \mathbb{N}\cup \{0\}$, the following formulas hold:
\begin{align}\label{c2kn_exact}
	c(m, 2j, n) &= \dfrac{2^{(n-1)/2}\pi\Gamma\left(
		j+\frac{n}{2}
	\right) \lambda_{j,n}}{\Gamma\left(
		j + \frac{1}{2}
	\right)}, \\ 
	\label{thm_even_const_lambda}
	\lambda_{j,n} &= 2(-1)^{j}
	\left[
		\sqrt{\pi}\dfrac{\Gamma(j + \frac{1}{2})}{\Gamma(j + 1)}
	\right]^{n-1},
\end{align}
where $c(m,k,n)$ are the constants arising in \eqref{integral_claims}, $\Gamma(\cdot)$ is the Gamma function, $\lambda_{j,n}$ are the eigenvalues of the Minkowski-Funk operator $\mathcal{M}$ defined in \eqref{minkowski_funk}.
\end{thm}
In the case of $n=1, \, k=2j$ formulas \eqref{integral_claims}-\eqref{integral_sign_claims} arise in formula (7) of \cite{kunyansky1992generalized}. In the case of $n=2$ formulas \eqref{integral_claims}-\eqref{integral_sign_claims} and constants $c(m,k,n)$ for general $k$ were given in formula (1) of \cite{neves2006analytical}. 
\par We didn't success to find in literature the explicit values for the eigenvalues $\lambda_{j,n}$ of operator $\mathcal{M}$ defined in \eqref{mkn_def} for $n > 2$.
\par In particular, formulas \eqref{integral_claims}-\eqref{thm_even_const_lambda} are very essential for inversion of weighted Radon transforms; see \cite{kunyansky1992generalized} for $n=1$ and the subsequent work \cite{gon2016nov} for $n\geq 2$.
\begin{rem}\label{lem_general_sum}
Let $c(m,k,n)$ be the constant arising in \eqref{integral_claims} and function $\Phi(k,n)$ be defined as follows:
\begin{equation}\label{Phi_func_def}
	\Phi(k,n) \eqdef \sum\limits_{m=1}^{a_{k,n+1}} c(m,k,n),
\end{equation}
where $a_{k,n+1}$ is defined in \eqref{eigensp_dim}. Then function $\Phi(k,n)$ does not depend on the particular basis of spherical harmonics $\{Y_k^m\}$ for fixed $k,n$.
\end{rem}

\section{Proofs}\label{sect_proofs}
\subsection{Proof of formula \eqref{integral_claims}}
From the Funk-Hecke Theorem (see e.g. \cite{morimoto1998analytic}, Chapter 2, Theorem 2.39) it follows that:
\begin{equation}\label{int_funk_hecke}
	\int\limits_{\Sp^{n}} Y_{k}^{m}(\theta)e^{i\rho(p\theta)}\, d\theta = 
	Y_{k}^{m}(p)c_{k}^{m}(\rho), \, m = \overline{1,a_{k,n+1}}.
\end{equation}
Formula \eqref{integral_claims} follows from \eqref{int_funk_hecke}, the following  differential equation:
\begin{align}\label{cki_diff_def}
	&\dfrac{1}{\rho^{n}}\dfrac{d}{d \rho}
		\left(
			\rho^{n} \dfrac{d c_{k}^{m}}{d \rho}
		\right)
		+ \left(1-\dfrac{\mu_{k,n}}{\rho^2}\right) c_{k}^{m} = 0, \, 
		c_k^m(0) = 0 \text{ for } k \geq 1,\\ \label{mkn_def}
	&\mu_{k,n} = k(n+k-1), \, 	
\end{align}
where function $c_k^m(\rho)$ arises in the right hand-side of \eqref{int_funk_hecke}, 
and from the fact that the solution of equation \eqref{cki_diff_def} with indicated boundary condition is given by the formula:
\begin{equation}\label{solution_general_diff}
	c_m^k(\rho) = c(m,k,n) \rho^{(1-n)/2}J_{k + \frac{n-1}{2}}(\rho), 
\end{equation}
where $c(m,k,n)$ is some constant depending on integers $m,k,n$; $J_r(\rho)$ is the $r$-th Bessel function of the first kind (see e.g. \cite{temme2011special}).
\par Formulas \eqref{int_funk_hecke}, \eqref{solution_general_diff} imply formula \eqref{integral_claims}.
\par It remains to prove that formulas \eqref{cki_diff_def}-\eqref{solution_general_diff} hold.
First, we prove that formulas \eqref{cki_diff_def}, \eqref{mkn_def} hold.

\par We recall that Laplacian $\Delta$ in $\R^{n+1}$ in spherical coordinates is given by the formula:
\begin{equation}\label{laplace_def}
	\lap u = \dfrac{1}{\rho^n}\left(
		\rho^n \dfrac{d}{d\rho}u
	\right) + \dfrac{1}{\rho^2}\lap_{\Sp^n}u, \, \rho \in (0,+\infty),
\end{equation}
where $u$ is a test function.
\par From formulas \eqref{integral_def}, \eqref{eigen_sp_lap}, \eqref{int_funk_hecke}, \eqref{laplace_def} it follows that 
\begin{equation}\label{laplace_integral_1}
 \lap I_{k}^{m} (p,\rho) = Y_{k}^{m}(p)\left(
		\dfrac{1}{\rho^{n}}\dfrac{d}{d \rho}
		\left(
			\rho^{n} \dfrac{d c_{k}^{m}}{d \rho}
		\right)
		- \dfrac{\mu_{k,n}c_{k}^{m}}{\rho^2}
	\right),
\end{equation}
where $\mu_{k,n}$ is defined in \eqref{mkn_def}.
\par On the other hand,
\begin{align}\nonumber 
	\lap I_{k}^{m}(p,\rho) &= \int\limits_{\Sp^{n}} Y_{k}^{m}(\theta)
	\lap e^{i\rho(p\theta)}\, d\theta \\ \label{laplace_integral_2}
	&= - \int\limits_{\Sp^{n}} Y_{k}^{m}(\theta)\, |\theta|^2\, e^{i\rho (p\theta)}\, d\theta = - I_{k}^{m}(p,\rho), \, p\in \Sp^n, \, \rho \geq 0.
\end{align}
Formulas \eqref{cki_diff_def}, \eqref{mkn_def} follow from \eqref{integral_def}, \eqref{laplace_integral_1}, \eqref{laplace_integral_2}. In particular, the boundary condition in \eqref{cki_diff_def} follows from orthogonality of $\{Y_k^m\}$ on $\Sp^n$ and the following formulas:
\begin{align}
	&I_k^m(0,p) = \int\limits_{\Sp^n} Y_k^m(\theta) \, d\theta = \begin{cases}
		0, \, k \geq 1, \\
		\mathrm{vol}(\Sp^n)c, \, k = 0
	\end{cases},\\
	&Y_0^1(p) = c \neq 0, \, p \in \Sp^n,
\end{align}
where $I_k^m$ is defined in \eqref{integral_def}, $\mathrm{vol}(\Sp^n)$ denotes the standard Euclidean volume of $\Sp^n$.
\par Next, formula \eqref{solution_general_diff} is proved as follows.
\par We use the following notation for fixed $k,m$:
\begin{align}\label{y_notation}
	y(t) = c_{k}^{m}(\rho), t = \rho \geq 0.
\end{align}
Using differential equation \eqref{cki_diff_def} and the notations from \eqref{y_notation} we obtain:
\begin{equation}\label{target_diff}
	t y''(t) + ny'(t) + \left(t-\dfrac{\mu_{k,n}}{t}\right) y(t) = 0, \, y(0) = 0,\, 
	\text{ for } k\geq 1.
\end{equation}
\par In order to solve \eqref{target_diff} we make the following change of variables:
\begin{equation}\label{chng_var}
	y(t) = t^{(1-n)/2}Z(t), \, t\geq 0.
\end{equation}
Formula \eqref{chng_var} implies the following expressions for $y'(t), \, y''(t)$ arising in \eqref{target_diff}:
\begin{align}\label{y1_z}
	y'(t) &= \dfrac{1-n}{2}t^{-(1+n)/2}Z(t) + t^{(1-n)/2}Z'(t),\\ 
	\label{y2_z}
	y''(t) &= \dfrac{(n^2-1)}{4}t^{-(1+n)/2}\dfrac{Z(t)}{t} + 
	(1-n)t^{-(1+n)/2}Z'(t) + t^{(1-n)/2}Z''(t), \, t\geq 0,
\end{align}
where $Z(t)$ is defined in \eqref{chng_var}.
\par Using formulas  \eqref{target_diff}, \eqref{y1_z}, \eqref{y2_z} we obtain:
\begin{equation}\label{z_diff_eq}
	t Z''(t) + Z'(t) + 
	\left(
		t - \dfrac{\left(k + \frac{n-1}{2}\right)^2}{t}
	\right)Z(t) = 0, \, t\geq 0.
\end{equation}
Differential equation \eqref{z_diff_eq} for unknown function $Z(t)$ is known as Bessel differential equation of the first kind with parameter $k + (n-1)/2\in \R$ (see e.g. \cite{temme2011special}). The complete solution of \eqref{z_diff_eq} is given by the following formula:
\begin{equation}
	Z(t) = C_1 J_{k + \frac{n-1}{2}}(t) + C_2 Y_{k + \frac{n-1}{2}}(t),\, t\geq 0,
\end{equation}
where $J_{r}(t), Y_r(t)$ are $r$-th Bessel functions of the first and second kind, respectively, $C_1, C_2$ are some constants; see e.g. \cite{temme2011special}.
Boundary condition in \eqref{target_diff} implies that 
\begin{equation}\label{z_sol}
	Z(t) = C_1 J_{k + \frac{n-1}{2}}(t), \, t\geq 0.
\end{equation}
\par Formulas \eqref{y_notation}, \eqref{chng_var}, \eqref{z_sol} imply that \eqref{solution_general_diff} is the complete solution of \eqref{target_diff}.
\par Formula \eqref{integral_claims} is proved.

\subsection{Proof of formula \eqref{integral_sign_claims}}
\par Formula \eqref{integral_sign_claims} follows from definition \eqref{integral_def} and the following property of the spherical harmonic $Y_k^m$:
\begin{equation}\label{harmonic_symm}
	Y_k^m(-\theta) = (-1)^kY_k^m(\theta), \, \theta\in \Sp^n.
\end{equation}
Property \eqref{harmonic_symm} reflects the fact that $Y_k^m(\theta) = Y_k^m(\theta_1, \dots, \theta_{n+1})$ is a homogeneous polynomial of degree $k$ restricted to $\Sp^n$. 
\par Formula \eqref{integral_sign_claims} in Theorem~\ref{thm_target} is proved.

\subsection{Proof of formulas \eqref{c2kn_exact}, \eqref{thm_even_const_lambda}}
	Formula \eqref{c2kn_exact} follows from orthonormality of $\{Y_k^m\}$ (in the sense of $L^2(\Sp^n)$), from formulas \eqref{integral_def}, \eqref{integral_claims}, \eqref{integral_sign_claims}, the following formula:
	\begin{equation}
		\label{eigen_funk_prop}
		\mathcal{M}[Y_{2k}^m] = \lambda_{k,n} Y_{2k}^m
	\end{equation}
	and from the following identities:
	\begin{align}\label{proof_even_const_1}
		\int\limits_{\Sp^n}\overline{Y_{2k}^m}(p) \, dp
		\int\limits_{-\infty}^{+\infty} I_{2k}^m(p, \rho) \, d\rho &= 2c(m, 2k, n) 
		\int\limits_{\Sp^n} |Y_{2k}^m(p)|^2 dp
		\int\limits_{0}^{+\infty} J_{2k + \frac{n-1}{2}}\rho^{(1-n)/2} d\rho \\ \nonumber
		&= 2c(m,2k,n) \int\limits_{0}^{+\infty} J_{2k + \frac{n-1}{2}}\rho^{(1-n)/2} d\rho \\ \nonumber
		&= c(m,2k,n) \dfrac{2^{(3-n)/2}\Gamma(\frac{1}{2}+k)}{\Gamma(k + \frac{n}{2})},\\
		\label{proof_even_const_2}
		\int\limits_{\Sp^n}\overline{Y_{2k}^m}(p)
		\int\limits_{-\infty}^{+\infty} I_{2k}^m(p, \rho) \, d\rho &= 
		2\pi
		\int\limits_{\Sp^n}\overline{Y_{2k}^m}(p)\, dp\int\limits_{\Sp^n}Y_{2k}^m(\theta)\delta(p\theta)\,d\theta \\ \nonumber
		&=2\pi \int\limits_{\Sp^n} \overline{Y_{2k}^m}(p) \mathcal{M}[Y_{2k}^m](p)\, dp = 
		2\pi\lambda_{k,n} \int\limits_{\Sp^n}|Y_{2k}^m|^2(p)\,dp\\ \nonumber
		&= 2\pi\lambda_{k,n}, 
\end{align}	
where $c(m,2k,n)$ arises in \eqref{integral_claims}, $\delta = \delta(s)$ is 1D Dirac delta function, $\overline{Y_k^m}$ is the complex conjugate of $Y_k^m$, $\mathcal{M}[Y^m_{k}]$ is  defined in \eqref{minkowski_funk}, $\lambda_{k,n}$ is given in \eqref{thm_even_const_lambda}, 
$\Gamma(\cdot)$ is the Gamma function.
	\par Formula \eqref{eigen_funk_prop} reflects the known property of the Funk-Minkowski transform $\mathcal{M}$ that the eigenvalue $\lambda_{k,n}$ of operator $\mathcal{M}[\cdot]$ defined in \eqref{minkowski_funk} corresponds to the eigensubspace of harmonic polynomials of degree $2k$ on $\R^{n+1}$ restricted to $\Sp^n$ (see e.g. \cite{gel_fand2003selected}, Chapter 6, p. 24). Note that, in \cite{gel_fand2003selected} it was proved that formula \eqref{eigen_funk_prop} holds also for all harmonic polynomials in $\R^3$ (i.e. $n=2$) restricted to $\Sp^2$, however these considerations admit a straightforward generalization to the case of arbitrary dimension $n \geq 1$.
	\par Formulas \eqref{proof_even_const_1}, \eqref{proof_even_const_2} imply formula \eqref{c2kn_exact}.
	\par Now, it remains to find the explicit value for $\lambda_{k,n}$ in formula \eqref{eigen_funk_prop}. We obtain it according to \cite{gel_fand2003selected} (Chapter 2, page 24), where the case of dimension $n=2$ was considered. 
	\par In particular, formula \eqref{eigen_funk_prop} holds for any homogeneous harmonic polynomial $P_{2k}$ of degree $2k$ in $\R^{n+1}$, where $P_{2k}$ is restricted to the sphere $\Sp^n \subset \R^{n+1}$; see \cite{gel_fand2003selected} (Chapter 2).
	\par We consider the following harmonic polynomial in $\R^{n+1}$:
	\begin{equation}\label{harmonic_polynom}
		P_{2k}(x) = P_{2k}(x_1, \dots, x_{n+1}) = (x_{n} + ix_{n+1})^{2k}, \, x \in \R^{n+1}.
	\end{equation}
	From formula \eqref{harmonic_polynom}, aforementioned results of \cite{gel_fand2003selected} (Chapter 6) and their generalizations to the case of arbitrary dimension $n \geq 1$ it follows that $P_{2k}$ being restricted to sphere $\Sp^n$ is an eigenfunction of operator $\mathcal{M}$ which corresponds to the eigenvalue $\lambda_{k,n}$. 
	\par We consider the spherical coordinates in $\R^{n+1}$ given by the following formulas:
	\begin{align}\nonumber
		x_1 &= \cos(\theta_n),\\ \nonumber
		x_2 &= \sin(\theta_n)\cos(\theta_{n-1}),\\ \label{sph_coord}
		\cdots\\ \nonumber
		x_{n} &= \sin(\theta_n) \sin(\theta_{n-1})\cdots \sin(\theta_2)\cos(\phi),\\ \nonumber
		x_{n+1} &= \sin(\theta_n) \sin(\theta_{n-1})\cdots \sin(\theta_2)\sin(\phi),\\ \nonumber
		\theta_n&, \, \theta_{n-1}, \cdots, \theta_2 \in [0,\pi], \, \phi \in [0,2\pi).
	\end{align}
	Formulas  \eqref{harmonic_polynom}, \eqref{sph_coord} imply that polynomial $P_{2k}$ being  restricted to $\Sp^n$ may be rewritten as follows:
\begin{equation}\label{spher_polynom_angle}
	P_{2k}|_{\Sp^n} = P_{2k}(\theta_n, \, \theta_{n-1}, \dots, \theta_1, \phi) = e^{i2k\phi} \prod_{i=2}^{n} \sin^{2k}(\theta_i),
\end{equation}
where $(\theta_n, \dots, \theta_1, \phi)$ are the coordinates on $\Sp^n$ according to \eqref{sph_coord}.

	From formulas in \eqref{harmonic_polynom}, \eqref{sph_coord}, \eqref{spher_polynom_angle} it follows that:
	\begin{equation}\label{spher_polynom_atpoint}
		P_{2k}|_{\Sp^n}= P_{2k}(\pi/2, \dots, \pi/2, 0) = 1.
	\end{equation}
\par From formulas \eqref{eigen_funk_prop}, \eqref{spher_polynom_atpoint}, \eqref{spher_polynom_angle} we obtain:
	\begin{align} \nonumber
		\lambda_{k,n} = \mathcal{M}[P_{2k}](\pi/2, \dots , \pi/2, 0) &= 
		2(-1)^{k}\prod_{i=2}^{n}\int\limits_{0}^{\pi} \sin^{2k}(\theta_i) \, d\theta_i \\
		&= 2(-1)^k
		\left[
			\sqrt{\pi}\dfrac{\Gamma(k + \frac{1}{2})}{\Gamma(k + 1)}
		\right]^{n-1}, 
	\end{align}
	which implies \eqref{thm_even_const_lambda}.
	\par Formulas \eqref{c2kn_exact}, \eqref{thm_even_const_lambda} are proved.
	
\subsection{Proof of Remark~\ref{lem_general_sum}}
The statement of the remark follows from the following formulas:
\begin{align}
	\label{lem_doubl_integrand}
	\sum\limits_{m=1}^{a_{k,n+1}} c(m,k,n) &= 
	\left(
		\rho^{(1-n)/2}J_{k + \frac{n-1}{2}}(\rho)
	\right)^{-1} 
	\dfrac{a_{k,n+1}}{\vol(\Sp^n)}
	\int\limits_{\Sp^n}\int\limits_{\Sp^n}P_{k,n}(p\theta)e^{i\rho(p\theta)}\, dp\,d\theta,\\
	\rho \geq 0 &\text{ such that }\rho^{(1-n)/2}J_{k + \frac{n-1}{2}}(\rho) \neq 0, \\
	\label{gegen_inv_def}
	P_{k,n}(\sigma\theta) &= \dfrac{\vol(\Sp^n)}{a_{k,n+1}}\sum\limits_{m=1}^{a_{k,n+1}}
	Y_k^m(\sigma)\overline{Y_{k}^m(\theta)}, \, \sigma,\theta\in \Sp^n,
\end{align}
where $a_{k,n+1}$ is defined in \eqref{eigensp_dim}, $\mathrm{vol}(\Sp^n)$ denotes the standard Euclidean volume of $\Sp^n$, $\overline{Y_k^m}$ is the complex conjugate of $Y_k^m$, $P_{k,n}(t), \, t\in [-1,1]$ is a Gegenbauer's polynomial defined by the formula (see \cite{morimoto1998analytic}):
\begin{equation}\label{pol_gegen_def}
	P_{k,n}(t) \eqdef \dfrac{(-1)^k}{2^k}
	\dfrac{\Gamma(\frac{n}{2})}{\Gamma(k+\frac{n}{2})}
	\dfrac{1}{(1-t^2)^{\frac{n-2}{2}}}
	\dfrac{d^k}{dt^k}(1-t^2)^{k + \frac{n-2}{2}}, \, k\in \mathbb{N}\cup\{0\}, \, n\geq 2,
\end{equation}
where $\Gamma(\cdot)$ is the Gamma function.
\par The integrand in $dp\, d\theta$ in the right hand-side of \eqref{lem_doubl_integrand} is independent of basis $\{Y_{k}^m\}$, which implies the statement of Remark~\ref{lem_general_sum}.
\par It remains to prove that formulas \eqref{lem_doubl_integrand}, \eqref{gegen_inv_def} hold,  where $P_{k,n}(t), \, t\in [-1,1]$ is defined in \eqref{pol_gegen_def}.
\par Formula \eqref{gegen_inv_def} is proved in \cite{morimoto1998analytic} (Chapter 2),
which is also the known property of Gegenbauer's polynomials.
\par From formulas \eqref{integral_def}, \eqref{integral_claims}, \eqref{gegen_inv_def} it follows that:
\begin{align}\nonumber \label{lem_gen_sum_preint}
	\dfrac{\vol(\Sp^n)}{a_{k,n+1}}\sum \limits_{m=1}^{a_{k,n+1}} \int\limits_{\Sp^n}I_{k}^m(p,\rho) \overline{Y_{k}^m(p)} \, dp &=
	\int\limits_{\Sp^n} 
	\int\limits_{\Sp^n} P_{k,n}(p \theta) e^{i\rho(p\theta)} d\theta\, dp \\ \nonumber
	&= \rho^{(1-n)/2}J_{k + \frac{n-1}{2}}(\rho)\dfrac{\vol(\Sp^n)}{a_{k,n+1}} \sum\limits_{m=1}^{a_{k,n+1}}
	c(m,k,n) \int\limits_{\Sp^n}|Y_{k}^m(p)|^2 dp,\\
	&= \rho^{(1-n)/2}J_{k + \frac{n-1}{2}}(\rho)\dfrac{\vol(\Sp^n)}{a_{k,n+1}} \sum\limits_{m=1}^{a_{k,n+1}}
	c(m,k,n).
\end{align}
where $\{Y_k^m\}$ is the orthonormal basis of spherical harmonics on $\Sp^n$, $P_{k,n}(t), \, t\in [-1,1]$ is defined in \eqref{pol_gegen_def}.
\par Formulas \eqref{gegen_inv_def}, \eqref{lem_gen_sum_preint} imply formula \eqref{lem_doubl_integrand}.
\par Remark~\ref{lem_general_sum} is proved.

\section{Aknowledgments}
The present work was fulfilled in the framework of research conducted under the direction of Prof. R. G. Novikov.
\par The author would like also to thank assistant-professor Hendrik De Bie at Ghent University who noticed the first version of this work and provided useful links to his own result on the case odd $k$ in Theorem 1; see \cite{hendr2007super}.


\begin{thebibliography}{NPF{\etalchar{+}}06}

\bibitem[A13]{adkins2013three}
G. S. Adkins.
\newblock { \em Three-dimensional Fourier transforms, integrals of spherical Bessel functions, and novel delta function identities}.
\newblock arXiv preprint arXiv:1302.1830, 2013.

\bibitem[AWH11]{arfken2011mathematical}
G. B. Arfken, H. J. Weber, F. E. Harris.
\newblock {\em Mathematical methods for physicists: a comprehensive guide},
\newblock Academic Press, 2011.


\bibitem[Fu13]{funk1913}
P. Funk.
\newblock {{\"U}ber Fl{\"a}chen mit lauter geschlossenen geod{\"a}tischen Linien},
\newblock {\em Mathematische Annalen}, 74(2):278-300, Springer,
\newblock {1913}

\bibitem[Gon16]{gon2016nov}
F. O. Goncharov.
\newblock {\em Iterative inversion of weighted Radon transforms in 3D},
\newblock preprint, 2016.


\bibitem[GGG03]{gel_fand2003selected}
I. M. Gel'fand, S. G. Gindikin, M. I. Graev.
\newblock {\em Selected topics in integral geometry}, vol. 220.
\newblock American Mathematical Society, 2003.

\bibitem[GN16]{goncharov2016analog}
F. O. Goncharov, R. G.~Novikov.
\newblock An analog of {C}hang inversion formula for weighted {R}adon
  transforms in multidimensions.
\newblock {\em Eurasian Journal of Mathematical and Computer Applications}, 4(2):23-32, 2016.


\bibitem[Kun92]{kunyansky1992generalized}
L. A. Kunyansky.
\newblock {\em Generalized and attenuated {R}adon transforms: restorative approach
  to the numerical inversion}.
\newblock {Inverse Problems}, 8(5):809, 1992.

\bibitem[Mor98]{morimoto1998analytic}
M. Morimoto.
\newblock {\em Analytic functionals on the sphere}.
\newblock American Mathematical Society, 1998.

\bibitem[NPF{\etalchar{+}}06]{neves2006analytical}
A. A. Neves, L.A. Padilha, A. F., E.
  Rodriguez, C. Cruz, L. Barbosa, C. Cesar.
\newblock Analytical results for a {B}essel function times Legendre polynomials
  class integrals.
\newblock {\em Journal of Physics A: Mathematical and General}, 39(18):293,
  2006.

\bibitem[SW16]{stein2016introduction}
E. M. Stein, G. Weiss.
\newblock {\em Introduction to {F}ourier analysis on {E}uclidean spaces
  (PMS-32)}, vol.~32.
\newblock Princeton University Press, 2016.

\bibitem[Tem11]{temme2011special}
N.~M. Temme.
\newblock {\em Special functions: {A}n introduction to the classical functions
  of mathematical physics}.
\newblock John Wiley \& Sons, 2011.

\bibitem[Hd07]{hendr2007super}
H. De Bie, F. Sommen.
\newblock Spherical harmonics and integration in superspace.
\newblock{\em Journal of Physics A: Mathematical and Theoretical}, 
26(40):7193, IOP Publishing, 2007.


\end{thebibliography}
\newcommand{\etalchar}[1]{$^{#1}$}

\end{document}